# Modeling the Complexity of City Logistics Systems for Sustainability

Taiwo Adetiloye[#1], Anjali Awasthi[#2]
[#]Department of Information and Systems Engineering
Concordia University
Montreal, Quebec, Canada
[1]t_adeti@encs.concordia.ca
[2]anjali.awasthi@concordia.ca

*Abstract*— The logistics of urban areas are becoming more sophisticated due to the fast city population growth. The stakeholders are faced with the challenges of the dynamic complexity of city logistics(CL) systems characterized by uncertainty effect together with the freight vehicle emissions causing pollution. In this conceptual paper, we present a research methodology for the environmental sustainability of CL systems that can be attained by effective stakeholders' collaboration under non-chaotic situations and the presumption of the human levity tendency. We propose the mathematical axioms of the uncertainty effect while putting forward the notion of condition effectors, and how to assign hypothetical values to them. Finally, we employ a spider network and causal loop diagram to investigate the system's elements and their behavior over time.

## I. INTRODUCTION

Rural-to-urban migration is rapidly increasing with most of the human populations living in the cities, today. It is expected that global urbanization would significantly increase by over 3 billion people in the next 30years [1]. Thus, the participants in city logistics (CL) operations are drawn towards adopting new approaches to overcome the traffic congestion caused by freight vehicle movements, to improve vehicle utilization, and to reduce emissions and pollutions without penalizing the city social and economic activities [2]. A number of research innovations have attempted to address the problems of CL for environment sustainability [3-6]. Some of these designs seek to fully optimize the logistics and transport activities by private companies in urban areas while considering the traffic environment, the traffic congestion, and energy consumption within the framework of a market economy [7]. Central to the issue at hand is the environmental concerns that need conscious and respectable attitude by the stakeholders with their strategies tailored towards green logistics [8]. In addition, efforts at improving goods transportation in a city have to take into account the dynamic complexity of a city and the activities of the stakeholders. Ref. [9] observed that despite the evolution of internet e-collaboration tools that play the vital role of "value creation enabler", companies still have to face the business and market "complexities".

It can be assumed that the change in the CL systems complexity are widely determined by variable factors that include the administrative policies, the activities of shippers, the information technology, the infrastructures in place, the residents' socio-cultural characteristics with demands, the freights, the goods, and the environment; while, information sharing and consolidations prioritize design innovations of CL given the backdrop of its dynamic complexity. The main objective of this paper is to investigate the complexity of CL systems from uncertainty perspective and to model its effect on the activities of the CL stakeholders and the city environment; the relevance of which can be disruptive to the supply chain if not alleviated through effective collaboration strategies.

### A. Problem definition

The prospects of designing a perfect CL systems given the complexity of its supply chain may be hindered by uncertainty. We assume that the CL systems complexity will either increase or decrease due to redundancy or emergence of c-commerce subsystems with respect to the change in the complexity of the city. From the foregoing, we attempt to answer the following:
- What critical assumptions are important for studying the complexity of CL systems?
- How can the uncertainty elements of the CL systems be used for modeling complexity?
- Can collaboration mitigate the uncertainty inherent in the complexity of CL systems?

### B. Practical significance

In real-world, simplifying the complexity model becomes necessary while making the presumption of human levity tendency (HLT) and assuming that the city is in a non-chaotic state. The macro or strategic level offers various possibilities for collaborations by the stakeholders which can be considered with a view to improve the CL activities through consolidation centres and social network. The key strategic decision-making process become useful at the micro-level of the CL system where the operation participants have to consider their vehicle selection, goods to vehicle assignment, goods distribution, and environmental impact assessment for efficient and optimal operations. Thus, the stakeholders are better positioned to seek viable business strategies within a collaborative framework in order to mitigate the uncertainty effect inherent in the complexity of CL systems [5, 6]. Also, evaluating partners for CL collaboration under sustainable freight regulations using BOCR and fuzzy GRA to generate rankings for collaboration partners has been proposed by [10]. This could be expanded to explain the uncertainty elements that can have positive or negative impacts on such collaborative efforts.

## II. Literature review

Transportation of goods in urban areas constitutes a major enabling factor to the economic and social activities of urban lives [11]. The daily services provided by the shippers, freight carriers, and administrators contribute to the complexity of the CL framework within the greater supply chain, bundled with the uncertainty that emanates in various forms [12]. The term 'complexity' tends to imply something with different parts in intricate arrangement [13]. There is the notion that complexity can be disorganized or organized with respect to systems. A disorganized complexity has problematic states of randomness present in a large number of variables in which each of the many variables exhibits individual erratic behaviors, but the system itself possesses a degree of orderliness or analyzable average properties. On the other hand, organized complexities are all problems that involve dealing simultaneously with a sizable number of factors which have systemic relationship [14]. Improper handling of complexity across the CL systems can lead to increased uncertainty, risk and consequently unanticipated cost [15].

Ref. [15] used cluster analysis and quantitative arguments to conclude that the traffic volume of the road, carbon dioxide emissions from road traffic, and Gross Domestic Product (GDP) are factors that could be used for dematerialization and materialization of a country. Ref. [2] considered [17, 18] in modeling the demand uncertainty in two tiered CL (2T-CL) systems (tactical) planning. In the 2T-CL systems, different vehicle fleets from two layers of facilities move loads among facilities and from them to appropriate customers. Their aim was to address complex traffic cases by providing timely delivery of loads to customer in an economically and environmentally efficient operations while utilizing urban vehicles and city freighters on the city streets and at satellites, as germane to the uncertainty and tactical planning. In addition, they offered the solutions of a two-stage stochastic model with four different recourse strategies and formulations that adapt the plan to the observe demand. Ref. [10] recommend the use of fuzzy BOCR-GRA approach for collaboration partner selection as an important tool to achieve collaboration planning among CL operators for operational efficiency under municipal freight regulations such as access, sizing and timing restrictions etc. They seek to minimize the complexity and uncertainty often associated with the CL system through effective collaboration. Ref. [12] proposed the main research directions for developing comprehensive methodology for supply chain network (SCN) design under uncertainty in the scope of SCN risk analysis, SCN hazards modeling, scenario development and sampling, value based SCN design models, modeling for robustness, and for resilience and responsiveness as well as solution methods. A modeling framework for a better understanding of complexity in supply chain has been presented by [19]. A review of the CL modelling efforts based on trends and gaps has been presented by [20]. It covers the complexity and diversity of current CL practice and the different aspects in the modeling selection process as seen from the perspectives of the stakeholders. Ref. [21] observe that system's complexity due to the long waiting time between human actions and environmental effect can be an extremely demanding challenge of governance, and often beyond the scope of a unilateral solution.

## III. Research methodology

### A. Complexity of CL systems

A widely accepted definition of CL by [7] is that it is "the process for totally optimizing the logistics and transport activities by private companies in urban areas while considering the traffic environment, the traffic congestion and energy consumption within the framework of a market economy". The complexities of CL system can be found in its various elements such as socio-cultural values of residents, administrative and government policies, and activities of shippers and freight carriers as well as the environment, infrastructures and information technologies. Also, market competitions among organizations cannot always be a viable option for optimum performance within the city and as such has to be minimized to a reasonable level.

It can be assumed that the individual and interrelated complexities of the specific aspects of CL(the administrative policies, the activities of shippers, the information technology, the infrastructures in place, the residents' socio-cultural characteristics and demands, the freights, the goods, and the environment) change with time under the subjective influence of uncertainty. Therefore, they will either increase or decrease due to emergence or redundancy of collaborative subsystems as the complexity of the city changes. This is given by:

$$S_C = f(C) \qquad (1)$$

With $C = P_o \cup S \cup I_t \cup I \cup R \cup F_e \cup G \cup E \cup ...$

In (1), $S_C$ represents the system complexity defined as a function of the city complexity(C) where C can be defined as the agglomeration or union of the individual complexities of the administrative policies ($P_o$), shippers' activities($S$), information technology ($I_t$), infrastructures (I), residents' socio-cultural characteristics with demands (R), freights ($F_e$), goods (G), and environment (E) etc. of a city. Also, it follows that:

$$\Delta S_C \propto \Delta C \qquad (2)$$
$$\propto \Delta(P_o \cup S \cup I_t \cup I \cup R \cup F_e \cup G \cup E \cup ...)$$
$$\propto \Delta(P_o \cup S \cup I_t \cup I \cup R \cup F_e \cup G \cup E \cup ...)$$
$$= k_o (P_o \cup S \cup I_t \cup I \cup R \cup F_e \cup G \cup E \cup ...)$$

Equation (2) denotes that the change in the CL systems complexity, $\Delta S_C$, may be directly proportional to the city complexity, $\Delta C$, where $k_o$ is is the uncertainty effect defined as:

$$k_o = \{ x \in R; 0 < x < 1 \} \qquad (3)$$

And, in (3), it can assumed that the neutrality "1" and the nullity "0" of $k_o$ hardly exist; where digits 1 and 0 are taken to be binary.

Let

$$\Delta A = \Delta P_o \cup \Delta S \cup \Delta I_t \cup \Delta I \cup \Delta R \cup \Delta F_e \cup \Delta G \cup ... \quad (4)$$

with exclusion of $\Delta E$

However, by default, and if based on the HLT, some presumptions may arise that:

$$k_o \cong 1, \quad \Delta A \cong 0 \quad and \quad \Delta E \cong 1 \quad (5)$$

The HLT is defined as the behavioral way of assuming all is normal within a human set boundary despite contrary evidence pointing to events, elsewhere, outside the boundary. Three conditions (or trio conditionality) emerge as follow:

*1) Lemma 1*
For a city in non-chaotic situation, that is in a normal state:

$$\frac{\Delta S_c}{k_o(\Delta A \cup \Delta E)} \cong 1 \quad (6)$$

A non-chaotic situation is a normal state that city dwellers expect to have minor accidents from natural and man-made disasters. In this situation, stakeholders through collaborative planning can achieve sustainable CL operations.

*2) Lemma 2*
For a city in a near chaotic situation (Non-approximate):

$$0 < \frac{\Delta S_c}{k_o(\Delta A \cup \Delta E)} < 1 \quad (7)$$

A near chaotic situation would arise from a natural disaster of huge proportion such as a hurricane, a tsunami, a high magnitude earthquake and so on, which can engulf a city for a short time lasting from few hours to few days.

*3) Lemma 3*
For a city in cataclysmic chaotic state:

$$\frac{\Delta S_c}{k_o(\Delta A \cup \Delta E)} \cong 0 \quad (8)$$

This can result from high catastrophe such as an atomic explosion, a planetoid, a black hole and so on, which would unleash uncontrollable chain reactions that can engulf a city and other cities far and near.

The approximate equalities of one and zero for the equations contained in (5), (6) and (8) give credence to the law of conservation of mass which states that the mass of an isolated system will remain constant over time. Otherwise, if we are to assume exact equalities of greater than one or an exact equality of zero for these equations then the law of conservation of mass will be nullified with the implication that a city can be totally sustained, that is created or destroyed.

A detailed practical application of this model can be found in [6] based on stakeholders' collaboration for sustainable CL operation involving vehicle selection, goods to vehicle assignment, goods distribution and environmental impact assessment at the micro or operational level of CL operations. It should be noted that while their collaboration square model simplified the complexity model presented in this work, their operational model extended the CILOSIM framework in [3]. Furthermore, the technique at the macro or strategic level of the CL operations is designed to evaluate the social-cultural characteristics, economy and environment impacts of the activities of city dwellers while ignoring the influence of uncertainty. For fuller discussion, the reader should please refer to [3] and [5].

*B. Uncertainty and CL systems complexity*

A mathematical diagnosis of $k_o$ revealed that it can be evaluated by assigning binary values to the unpredictability effect of the (natural) conditions such as air, wind, snow, and dry weather, as well as that of tornado, hurricane, planetoid (hitting a city), alien invasion and so on. The validity of this argument is based on two basic axioms:

- The neutrality "1" and nullity "0" of $k_o$ hardly exist:

$$k_o = \{x \in R; 0 < x < 1\} \quad (9)$$

- $k_o$ is an effector (or effect vector) such that:

$$k_o = \left|\vec{k_o}\right| = \left|\vec{k_1}n_1 + \vec{k_2}n_2 + \vec{k_3}n_3 + \vec{k_4}n_4 + ...\right| \quad (10)$$

We say that the effectors of $k_o$ is the summation of individual condition effectors, $k_1 n_1, k_2 n_2, k_3 n_3, k_4 n_4$ etc. It is important to point out that an effector or effect vector is different from a vector because a vector has magnitude and direction but an effector has magnitude and positive or negative effects. The constraints necessary for axiom 2 to satisfy axiom 1 are:

$$k_1, k_2, k_2, k_4, ... = \{x \in R; 0 < x < 1\} \quad (11)$$

Scalars $n_1, n_2, n_3, n_4...$; cannot all be equal to zero, simultaneously

(12)

An effector can be either positive or negative. A condition can be defined as matters within the earth sphere such as air, water, wind, nuclear energy etc. found within the earth sphere, perceived matters such as sunlight, planetoid, aliens and matters that exist beyond our earth sphere, somewhere, in outer space such as in black holes etc. Conditions are identified as both negative and positive effectors. Conditions such as tornado, hurricane, planetoid, and black hole are primarily negative effectors. Hence, it can be stated that a negative effector causes a harmful, uncertainty effect while a positive effector causes an innocuous, uncertainty effect with impacts on the collaboration of CL stakeholders within a city. For example, a mild wind may have a negative effect if there is air pollution and with multiplicity – that is, a strong wind – it can endanger a city, thereby militating against collaboration of stakeholders involved in CL operations. Positive effectors often do not endanger a city; rather, they facilitate collaboration of stakeholders involved in CL operations but can be inhibitive with very high multiplicity. The following TABLE I describes the "hypothetic values" assigned for some of the condition effectors. The hypothetic values assigned for the conditions are based on the following assumptions:

*1)* *The human perception of the severity of a condition. For instance, air is perceived as least severe –hence, it is assigned the smallest hypothetic value –while black hole is perceived as most severe with the highest hypothetic value.*
*2)* *Conditions deemed as opposite in nature, such as wet and dry weather, are assigned same hypothetic value.*
*3)* *Conditions on the same scale of effect such as light and planetoid which emits photons and heat energy are assigned the same binary numbers before the decimal point but are differentiated with different binary numbers after the decimal point.*
*4)* *The further a condition exist from the earth sphere the greater its hypothetical value.*

Also, binary numbers with "1s" as hypothetic values were chosen because the arithmetic of binary numbers are better solved by computers than by the human mental cognition system due to the boredom and difficulty that become evident as its arithmetic complications unfolds. Secondly, binary numbers "1" and "0" are used to depict switch states. Moreover, for the hypothetic values of condition effectors although initialize with "1s" as a representation of a unique perfect state, a mixture of "1s" and "0s" often emerge as complex number in the axioms of the uncertainty effect. In addition, attempting to solve an axiom of uncertainty problem with binary numbers tend toward an unpredictable complex numeric that best mimics the uncertainty effect. Through this discussion, we provide a sense of the effect of uncertainty on the individual complexities of a city, and how it influences the collaboration of stakeholders involved in CL operations.

TABLE I. HYPOTHETIC VALUES OF CONDITION EFFECTORS

| $x$ | Condition | Negative effector $-k_x n_x$ | Positive effector $k_x n_x$ |
|---|---|---|---|
| 1 | Air | $-1.1n_1$ | $1.1n_1$ |
| 2 | Dry | $-1.11n_2$ | $1.11n_2$ |
| 4 | Wet | $-1.11n_3$ | $1.11n_3$ |
| 5 | Wind | $-11.1n_4$ | $11.1n_4$ |
| 6 | Snow | $-111.1n_5$ | $111.1n_5$ |
| 7 | Water | $-111.11n_6$ | $111.11n_6$ |
| 8 | Tornado | $-1111.1n_7$ | – |
| 9 | Hurricane | $-1111.11n_8$ | – |
| . | . | . | . |
| . | . | . | . |
| . | . | . | . |
| $a-3$ | Planetoid | $-11111.11n_{a-3}$ | – |
| $a-2$ | Sunlight | $-11111.1111n_{a-2}$ | $11111.1111n_{a-2}$ |
| $a-1$ | Alien invasion | $-11111111.1n_{a-1}$ | $11111111.1n_{a-1}$ |
| $a$ | Black hole | $-11111111111.1n_a$ | – |

### C. Uncertainty in CL systems

A spider's web is a spiraling polygon for a good reason, which is to meet the needs of the spider; albeit, the dexterity at which the spider holds together its web (at the center) and yet preys an insect is unpredictable or uncertain. This statement provides an analogy for categorizing the elements of CL systems. Fig. 1 is a spider diagram showing the various uncertainties present in CL system. The uncertainty in CL system can arise from change in the individual complexities of a city ($\Delta C$) which includes:

1) *Change in the administrative and government policies($\Delta P_o$).*
2) *Change in shippers' activities($\Delta S$).*
3) *Change in information technology ($\Delta I_t$).*
4) *Change in infrastructures ($\Delta I$).*
5) *Change due to residents' social and cultural values with demands ($\Delta R$).*
6) *Change in freights($\Delta F_e$).*
7) *Change in goods ($\Delta G$).*
8) *Change in the environment ($\Delta E$) of a city.*

They can be represented by the nodes of the octagonal spider network. The bi-directional links between these complexities are classified as either tangible or intangible connectors. Tangible links are directly measurable while intangible links are not measurable since they lack lucid physical appearance. However both type of connectors follow a clockwise, anticlockwise or dual-wise direction that signifies the agglomeration of the individual dynamic complexities of the city. A node having two tangible paths directly at opposite sides can be called tangible node, which is typically identified as $\Delta R$.

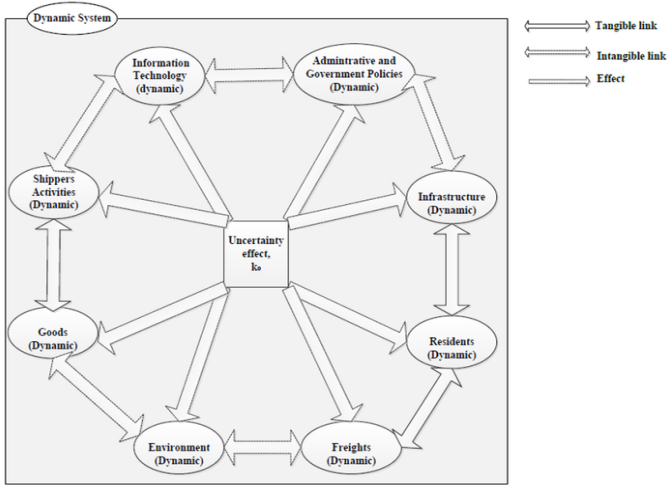

Fig. 1. Spider network showing linkage between uncertainty effect and individual complexities of a city

A node with two intangible paths directly at opposite sides can be seen as intangible, consisting $\Delta P_o$, $\Delta I_t$ and $\Delta E$ while a node with a tangible and an intangible link directly opposite sides can be called a semi-tangible node, which is the case for $\Delta S$, $\Delta G$, $\Delta I_t$ and $\Delta F_e$.

Unidirectional links extend from the black-box at the center towards each of the dynamic complexities located at the edges of the octagon. A unidirectional link is considered to be the influence that uncertainty effect exerts on the individual complexities of the city. The uncertainty effect represented by $k_o$ is depicted by the box at the center; while the dynamic system complexity, $\Delta Sc$, is depicted by the larger grey box enclosing the whole spider network. The components of the spider network is shown in TABLE II. Hence, the statement that the dynamic CL systems complexity, $\Delta S_C$, may be directly proportional to the dynamic complexity of the city, $\Delta C$.

TABLE II. COMPONENTS OF THE SPIDER NETWORK.

| Links | Node instance |
|---|---|
| Tangible | Residents |
| Intangible | Administrative and Government policies |
|  | Environment |
|  | Information Technology |
| Semi-tangible | Freights |
|  | Shippers |
|  | Goods |
|  | Infrastructure |

Furthermore, if we assume a dynamic CL systems state, $S_T$, defined as the mode or condition of a system, which is the reciprocal of the dynamic CL systems complexity as given by:

$$S_T = \frac{1}{\Delta S_C} \qquad (13)$$

Then, it can be deduced that:

$$\lim_{k_o \to \infty} S_T = 0 \qquad (14)$$

Equations (13) and (14) imply that a collaborative system, consisting of different subsystems such as B2B, B2G, C2B, G2B etc., ceases to exist if or whenever $k_o$ tends to infinity regardless of the agglomerative state of the individual complexities of the city.

### D. Collaboration and CL complexity

The strategies for collaboration rest on the need to create and expand the semi-intangible attributes and to optimize the use of intangible attributes of a city. This can be performed by fostering collaboration of stakeholders with efforts toward understanding the dynamic nature of individual city complexities such as the socio-cultural values of residents to demands, administrative and government policies, and activities of shippers and freight carriers as well as the environment, infrastructures and information technologies.

Evaluating these strategies might be based on how well the collaborative communities utilize the key knowledge of the axioms of uncertainty effect and the spider networks of CL. One approach for such evaluation is the use of causal loop diagram (CLD) for visualizing the relationship between the elements of a CL systems as illustrated in Fig. 2.

The link positive polarity (+) links points from one variable parameter to another and implied that an increase (or decrease) in input variable parameter would lead to an increase (or decrease) in the output variable parameter. For instance, increasing (or decreasing) changes in administrative policies as determined by tax rates, regulations and structuring directly lead to increasing (or decreasing) changes in shippers' activities and information technology. Likewise, links with positive polarities exist between shippers and freight carriers on the basis of supply rate of goods, and also between residents and information technologies. It can be observed that an increase in the activities of freight carriers and infrastructures contributes to increasing congestions within the city system. The term "congestion" in the context of CL can be defined as space minus freights and infrastructures. An important observation is that collaboration and pollutions have link positive polarities from nearly all the variables that extend from the administrators, freight carriers, shippers and residents.

Similarly, the link of negative polarity (-) implies that an increase (or decrease) in input variable parameter would lead to a decrease (increase) in the output variable parameter. The uncertainty effect projects a link negative polarity to almost all the variables identified in the CLD. The uncertainty effect is classified as mild and extreme. The difference between these two uncertainties has to do with the fact that there is always a delay associated with mild uncertainty effect while instantaneity is associated with extreme uncertainty effect. This means that increasing (or decreasing) mild and extreme uncertainty effect causes a decrease (or increase) in the dynamic states of the city complexities, which is delayed if the uncertainty effect is mild and instantaneous if the uncertainty effect is extreme. Simulating the realities of an uncertainty effect can be a daunting if not an impossible task.

Therefore, we can say that the interests of stakeholders are best assured by enhancing collaboration, as their numbers and activities increase in order for them to combat the growing challenges of environmental pollutions and the uncertainty effect.

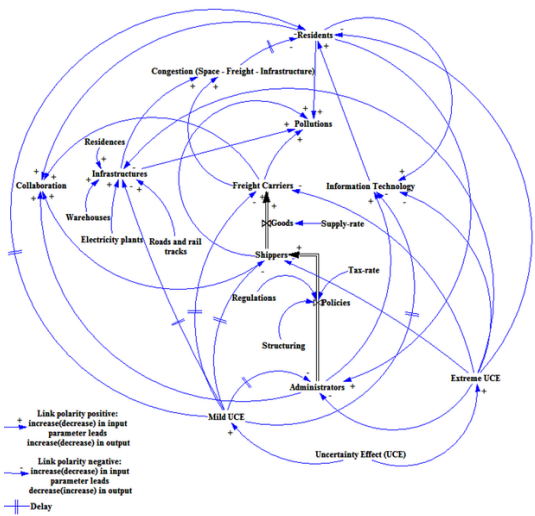

Fig. 2. Causal loop diagram for visualizing CL operations (designed with Vensim[1])

## IV. Conclusions

In this paper, we investigate the complexity of CL systems for achieving environmental sustainability. This understanding of CL complexity can assist decision makers in developing deeper insight into the impacts of supply chain activities on the city environment, emphasizing how the creation of efficient collaboration strategies alleviates uncertainty effect while ensuring that negative environmental impacts can be substantially reduced. Our future work involves detailed investigation of stakeholders' collaboration strategies for sustainability while respecting urban regulations and leveraging the advantages of ICT for decision making.

---

[1] Vensim is a registered trademark of the Ventana Systems, Inc. http://www.vensim.com/